\documentclass[12pt]{amsart}
\usepackage{verbatim}
\usepackage{amsmath,amssymb,amsbsy,amsfonts,amsthm,latexsym,
                        amsopn,amstext,amsxtra,euscript,amscd,mathrsfs,color,bm,cite}
\usepackage{mathtools}

\usepackage{todonotes}
\usepackage{url}
\usepackage[colorlinks,linkcolor=blue,anchorcolor=blue,citecolor=blue,backref=page]{hyperref}
\usepackage{todonotes}
\usepackage[norefs,nocites]{refcheck}

\usepackage{mathrsfs}

\usepackage{enumitem}

\newtheorem{theorem}{Theorem}[section]

\newtheorem{conjecture}[theorem]{Conjecture}

\def\({\left(}
\def\){\right)}
\def\[{\left[}
\def\]{\right]}
\def\<{\langle}
\def\>{\rangle}
\def\fl#1{\left\lfloor#1\right\rfloor}

\theoremstyle{plain}

\numberwithin{equation}{section}
\numberwithin{theorem}{section}
\numberwithin{table}{section}

\theoremstyle{remark}
\newtheorem{remark}[theorem]{Remark}

\newcommand{\Q}{{\mathbb Q}}

\newcommand{\Z}{{\mathbb Z}}

\newcommand{\OO}{\mathcal{O}}

\newcommand{\Qbar}{\bar{\Q}}

\newcommand{\bZ}{{\mathbb Z}}
\newcommand{\bG}{{\mathbb G}}
\newcommand{\bN}{{\mathbb N}}

\newcommand{\bF}{{\mathbb F}}
\newcommand{\Fq}{\bF_q}
\newcommand{\lra}{\longrightarrow}

\newcommand{\cS}{\mathcal{S}}

\title[Sparsity in the DML Conjecture]{A sparsity result for the Dynamical Mordell-Lang Conjecture in positive characteristic}

\author[D. Ghioca]{Dragos Ghioca}
\author[A. Ostafe]{Alina Ostafe}
\author[S. Saleh]{Sina Saleh}
\author[I. E. Shparlinski]{Igor E. Shparlinski}


\address{
Dragos Ghioca\\
Department of Mathematics\\
University of British Columbia\\
Vancouver, BC V6T 1Z2\\
Canada
}

\email{dghioca@math.ubc.ca}

\address{Alina Ostafe\\
School of Mathematics and Statistics\\
University of New South Wales\\
Sydney NSW 2052\\
Australia
}

\email{alina.ostafe@unsw.edu.au}

\address{
Sina Saleh\\
Department of Mathematics\\
University of British Co\-lumbia\\
Vancouver, BC V6T 1Z2\\
Canada
}

\email{sinas@math.ubc.ca}

\address{Igor E. Shparlinski\\ 
School of Mathematics and Statistics\\
University of New South Wales\\
Sydney NSW 2052\\
Australia
}

\email{igor.shparlinski@unsw.edu.au}

\begin{document}

\begin{abstract}
We prove a quantitative partial result in support of the Dynamical Mordell-Lang Conjecture (also known as the \emph{DML conjecture}) in positive characteristic. More precisely, we show the following: given a field $K$ of characteristic $p$,  given a semiabelian variety $X$ defined over a finite subfield of  $K$ and endowed with a regular self-map $\Phi:X\lra X$ defined over $K$, given a point $\alpha\in X(K)$ and a subvariety $V\subseteq X$, then the set of all non-negative integers $n$ such that $\Phi^n(\alpha)\in V(K)$ is a union of finitely many arithmetic progressions along with a subset $S$ with the property that there exists a positive real number $A$ (depending only on $N$, $\Phi$, $\alpha$, $V$) such that for each positive integer $M$, we have  
$$\#\left\{n\in S\colon~ n\le M\right\}\le A\cdot \left(1+\log M\right)^{\dim V}.$$
\end{abstract}

\maketitle


\section{Introduction}


\subsection{Notation}

Throughout this paper, we let $\bN_0:=\bN\cup\{0\}$ denote the set of nonnegative integers. 
As always in arithmetic dynamics, we denote by $\Phi^n$ the $n$-th iterate of the self-map $\Phi$ acting on some ambient variety $X$. For each point $x$ of $X$, we denote its orbit under $\Phi$ by 
$$\OO_\Phi(x):=\left\{\Phi^n(x)\colon~ n\in\bN_0\right\}.$$
Also, for us, an arithmetic progression is a set $\{an+b\}_{n\in \bN_0}$ for some $a,b\in\bN_0$; in particular, we allow the possibility that $a=0$, in which case, the above set is a singleton. 


\subsection{The Dynamical Mordell-Lang Conjecture}

The Dynamical Mordell-Lang Conjecture (see~\cite{GT-JNT}) predicts that for an endomorphism $\Phi$ of a quasiprojective variety $X$ defined over a field $K$ of characteristic $0$, given a point $\alpha\in X(K)$ and a subvariety $V\subseteq X$,  the set 
\begin{equation}
\label{eq:DML}
\mathcal{S}(\Phi, \alpha; V):=\left\{n\in\bN_0\colon~ \Phi^n(\alpha)\in V(K)\right\}
\end{equation} 
is a finite union of arithmetic progressions; for a comprehensive discussion of the Dynamical Mordell-Lang Conjecture, we refer the reader to the book~\cite{BGT}.

When the field $K$ has positive characteristic, then under the same setting as above, the return set $\mathcal{S}$ from~\eqref{eq:DML} is no longer a finite union of arithmetic progressions, as shown in~\cite[Examples~1.2 and~1.4]{G-TAMS}; instead, the following conjecture is expected to hold.

\begin{conjecture}[Dynamical Mordell-Lang Conjecture in positive characteristic]
\label{conj:DML}
Let $X$ be a quasiprojective variety defined over a field $K$ of characteristic $p$. Let $\alpha\in X(K)$, let $V\subseteq X$ be a subvariety defined over $K$, and let $\Phi:X\lra X$ be an endomorphism defined over $K$. Then the set $\mathcal{S}(\Phi, \alpha; V)$ given by~\eqref{eq:DML} is a union of finitely many arithmetic progressions along with finitely many sets of the form
\begin{equation}
\label{form of the solutions}
\left\{\sum_{j=1}^{m} c_j p^{a_j k_j}\colon~k_j\in\bN_0\text{ for each }j=1,\ldots m\right\},
\end{equation} 
for some given $m\in\bN$, some given $c_j\in\Q$, and some given $a_j\in\bN_0$ (note that in~\eqref{form of the solutions},  the parameters $c_j$ and  $a_j$ are fixed, while  
the unknowns $k_j$ vary over all non-negative integers, $j =1, \ldots, m$). 
\end{conjecture}

In~\cite{CGSZ}, Conjecture~\ref{conj:DML} is proven for regular self-maps $\Phi$ of tori assuming one of the following two hypotheses are met:
\begin{enumerate}
\item[(A)] $\dim V\le 2$;
\end{enumerate}
or
\begin{enumerate}
\item[(B)] $\Phi:\bG_m^N\lra \bG_m^N$ is a group endomorphism and there exists no nontrivial connected algebraic subgroup $G$ of $\bG_m^N$ such that an iterate of $\Phi$ induces an endomorphism of $G$ that equals a power of the usual Frobenius.
\end{enumerate}
The proof from~\cite{CGSZ} employs various techniques from Diophantine approximation (in characteristic $0$), to combinatorics over finite fields, to specific tools akin to semiabelian varieties defined over finite fields; in particular, the deep results of Moosa \& Scanlon~\cite{Moosa-Scanlon} are essential in the proof. 
Actually, the Dynamical Mordell-Lang Conjecture in positive characteristic turns out to be even more difficult than the classical Dynamical Mordell-Lang Conjecture since even the case of group endomorphisms of $\bG_m^N$ leads to deep Diophantine questions in \emph{characteristic $0$}, as shown in~\cite[Theorem~1.4]{CGSZ}. More precisely,~\cite[Theorem~1.4]{CGSZ} shows that solving Conjecture~\ref{conj:DML} just in the case of group endomorphisms of tori is \emph{equivalent} with solving the following polynomial-exponential equation: given any linear recurrence sequence $\{u_n\}$, given a power $q$ of the prime number $p$, and given positive integers $c_1,\ldots, c_m$ such that 
$$
\sum_{i=1}^m c_i<\frac{q}{2},
$$ 
then one needs to determine the set of all $n\in\bN_0$ for which we can find $k_1,\ldots ,k_m\in\bN_0$ such that
\begin{equation}
\label{eq:poly-expo 2}
u_n=\sum_{i=1}^m c_iq^{k_i}.
\end{equation}
The equation~\eqref{eq:poly-expo 2} remains unsolved for general sequences $\{u_n\}$ when $m>2$; for more details about these Diophantine problems, see~\cite{CZ3} and the references therein.


\subsection{Statement of our results}

Before stating our main result, we recall that a semiabelian variety is an extension of an abelian variety by an algebraic torus; for more details on semiabelian varieties, we refer the reader 
to~\cite[Section~2.1]{CGSZ} and the references therein. 

We prove the following result towards Conjecture~\ref{conj:DML}.
\begin{theorem}
\label{thm:main 2}
Let $K$ be a field of characteristic $p$, let $X$ be a semiabelian variety defined over a finite subfield of $K$, let $\Phi$ be a regular self-map of $X$ defined over $K$. Let $V\subseteq X$ be a subvariety defined over $K$ and let $\alpha\in X(K)$. Then 
the set $\cS(\Phi, \alpha; V)$ defined by~\eqref{eq:DML}  
is a union of finitely many arithmetic progressions along with a set $S\subseteq \bN_0$ for which there exists a  
constant $A$ depending only on $X$, $\Phi$, $\alpha$ and $V$ such that for all $M\in \bN$, we have
\begin{equation}
\label{eq:bound 3}
\#\left\{n\in S\colon~n\le M\right\}\le A\cdot \left(1 + \log M\right)^{\dim V}.
\end{equation} 
\end{theorem}

Our result strengthens~\cite[Corollary~1.5]{BGT-Funct} for the case of regular self-maps of semiabelian varieties defined over finite fields since  in~\cite{BGT-Funct} it is shown that the set $S$ (as in the conclusion of Theorem~\ref{thm:main 2}) is of Banach density zero; however, the methods from~\cite{BGT-Funct} cannot be used to obtain a sparseness result as the one from~\eqref{eq:bound 3}.  

We establish Theorem~\ref{thm:main 2} by combining~\cite[Theorem~3.2]{CGSZ} 
with~\cite[Th\'{e}or\`{e}me~6]{Laurent}.


\section{Proof of Theorem~\ref{thm:main 2}}

\subsection{Dynamical Mordell-Lang conjecture and linear recurrence sequences} 
First, since $X$ is defined over a finite field $\Fq$ of $q$ elements of characteristic $p$, we let $F:X\lra X$ be the Frobenius endomorphism corresponding to $\Fq$. We let $P\in \bZ[x]$ be the minimal polynomial with integer coefficients such that $P(F)=0$ in ${\rm End}(X)$; according to~\cite[Section~2.1]{CGSZ}, $P$ is a monic polynomial and it has simple roots $\lambda_1,\ldots, \lambda_\ell$, each one of them of absolute value equal to $q$ or $\sqrt{q}$. 

Using~\cite[Theorem~3.2]{CGSZ}, we obtain that the set $\cS(\Phi, \alpha; V)$ defined by~\eqref{eq:DML} is a finite union of \emph{generalized $F$-arithmetic sequences}, and furthermore, each such generalized $F$-arithmetic sequence is an intersection of finitely many \emph{$F$-arithmetic sequences}; see~\cite[Section~3]{CGSZ} for exact definitions.  
Each one of these $F$-arithmetic sequences consists of all non-negative integers $n$ belonging to a suitable arithmetic progression, for which there exist $k_1,\ldots,k_m\in\bN_0$ such that
\begin{equation}
\label{eq:F sequence 2}
u_n=\sum_{i=1}^m \sum_{j=1}^\ell c_{i,j}\lambda_j^{a_ik_i},
\end{equation}
for some given linear recurrence sequence $\{u_n\}\subset\Qbar$, some given $m,\ell\in\bN$, some given constants $c_{i,j}\in\Qbar$ and some given $a_1,\ldots, a_m\in\bN$. 
Applying Part~(1) of~\cite[Theorem~3.2]{CGSZ}, we also see that $m\le \dim V$. Furthermore, the linear recurrence sequence $\{u_n\}$ along with the constants $c_{i,j}$ and $a_j$ depend solely on $X$, $\Phi$, $\alpha$ and $V$.

Moreover, at the expense of further refining to another arithmetic progression, we may assume from now on, that the linear recurrence sequence $\{u_n\}$ is \emph{non-degenerate}, i.e the quotient of any two characteristic roots of this linear recurrence sequence is not a root of unity; furthermore, we may also assume that if one of the characteristic roots is a root of unity, then it actually equals $1$. For more details regarding linear recurrence sequences, we refer the reader to~\cite{Sch03}.   
In addition, we know that the characteristic roots of $\{u_n\}$ are all algebraic integers (see part~(2) of~\cite[Theorem~3.2]{CGSZ}); the characteristic roots of $\{u_n\}$ are either equal to $1$ (when $\Phi$ contains also a translation besides a group endomorphism) or equal to positive integer powers of the roots of the minimal polynomial of $\Phi$ inside ${\rm End}(X)$; for more details, see \cite[Section~3]{CGSZ}.  So, the equation~\eqref{eq:F sequence 2} becomes 
\begin{equation}
\label{eq:F sequence}
\sum_{r=1}^s Q_r(n)\mu_r^n = \sum_{i=1}^m\sum_{j=1}^\ell c_{i,j}\lambda_j^{a_ik_i},
\end{equation}
where $\mu_1,\ldots, \mu_s$ are the characteristic roots of the sequence $\{u_n\}$ and $Q_1,\ldots, Q_s\in \Qbar[x]$.

\subsection{Reduction to the case $s=1$} 
Now, if each polynomial $Q_r$ from the equation \eqref{eq:F sequence} is constant, then the famous result of Laurent \cite{Laurent} solving the classical Mordell-Lang conjecture (inside an algebraic torus) provides the desired conclusion that the set of all $n\in\bN_0$ satisfying an equation of the form~\eqref{eq:F sequence} must be a finite union of arithmetic progressions. So, from now on, we assume that not all of the polynomials $Q_r$ are constant.

Without loss of generality, we assume $Q_1$ is a non-constant polynomial. According to \cite[Section~8,~p.~319]{Laurent} (see also~\cite[Theorem~7.1]{Sch03}) all but finitely many solutions to the equation~\eqref{eq:F sequence} are also solutions to a \emph{subsum} corresponding to the equation~\eqref{eq:F sequence} which contains the term $Q_1(n)\mu_1^n$. More precisely, there exists a subset $1\in \Sigma_1\subseteq \{1,\dots, s\}$ and also, there exists a subset $\Sigma_2\subseteq \{1,\dots, m\}\times \{1,\dots, \ell\}$ such that
\begin{equation}
\label{eq:F sequence 3}
\sum_{r\in\Sigma_1} Q_r(n)\mu_r^n = \sum_{(i,j)\in \Sigma_2}c_{i,j}\lambda_j^{a_ik_i}.
\end{equation}
Moreover, letting $\pi_1:\{1,\dots, m\}\times \{1,\dots, \ell\}\lra \{1,\dots, m\}$ be the projection on the first coordinate, we have  $m_1:=\#\left(\pi_1(\Sigma_2)\right)$; in particular, $m_1\le m$. Without loss of generality, we assume $\pi_1(\Sigma_2)=\{1,\dots, m_1\}$ (with the understanding that, a priori, $m_1$ could be equal to $0$, even though we  show next that this is not the case).

Using \cite[Th\'{e}or\`{e}me~6]{Laurent}, the equation~\eqref{eq:F sequence 3} has finitely many solutions, unless the following subgroup $G_\Sigma\subseteq \Z^{1+m_1}$ is nontrivial. As described in \cite[Section~8,~p.~320]{Laurent}, the subgroup $G_\Sigma$ consists of all tuples $(f_0,f_1,\dots, f_{m_1})$ of integers with the property that 
\begin{equation}
\label{eq:same root}
\mu_r^{f_0}=\lambda_j^{a_if_i}\text{ for each $r\in\Sigma_1$ and each $(i,j)\in\Sigma_2$.}
\end{equation}  
Since $\mu_{r_2}/\mu_{r_1}$ is not a root of unity if $r_1\ne r_2$, we conclude that if $\Sigma_1$ contains at least  two elements (we already have by our assumption that $1\in\Sigma_1$), then $f_0=0$ in~\eqref{eq:same root}; furthermore, if $f_0=0$, then the equation~\eqref{eq:same root} yields  
that each $f_i=0$ (since each $\lambda_j$ has an absolute value greater than $1$ and $a_i\in\bN$). So, if $\Sigma_1$ has more than one element, then the subgroup $G_\Sigma$ is trivial and thus, \cite[Th\'{e}or\`{e}me~6]{Laurent} yields that the equation~\eqref{eq:F sequence 3} (and therefore, also the equation~\eqref{eq:F sequence}) has finitely many solutions, as desired.

\subsection{Concluding the argument} 
Therefore, from now on, we may assume that $\Sigma_1$ has a single element, i.e., $\Sigma_1=\{1\}$.  In particular, this also means that $\Sigma_2$ cannot be the empty set since otherwise the equation~\eqref{eq:F sequence 3} would simply read
$$Q_1(n)\mu_1^n=0,$$
which would only have finitely many solutions $n$ (since $\mu_1\ne 0$ and $Q_1$ is non-constant). So, we see that indeed $\Sigma_2$ is nonempty, which also means that $1\le m_1\le m$.

We have two cases: either $\mu_1$ equals $1$, or not.

{\bf Case 1.} $\mu_1=1$.

Then the equation~\eqref{eq:F sequence 3} reads:
\begin{equation}
\label{eq:F sequence 10}
Q_1(n)=\sum_{(i,j)\in\Sigma_2} c_{i,j}\lambda_j^{a_ik_i}.
\end{equation}
Now, for the equation~\eqref{eq:F sequence 10}, the subgroup $G_\Sigma$ defined above as in \cite[Section~8,~p.~320]{Laurent} is the subgroup $\Z\times \{(0,\dots, 0)\}\subset \Z^{1+m_1}$ since each integer $f_i$ from the equation~\eqref{eq:same root} must equal  $0$ for $i=1,\dots, m_1$ (note that $\mu_1=1$, while each $\lambda_j$ is not a root of unity). According to \cite[Th\'{e}or\`{e}me~6, part~(b)]{Laurent}, there exist positive constants $A_1$ and $A_2$ depending only on $Q_1$, the $c_{i,j}$ and the $a_i$ such that for any solution $(n,k_1,\dots, k_{m_1})$ of the equation~\eqref{eq:F sequence 10}, we have
\begin{equation}
\label{eq:logarithmic growth}
\max\left\{|k_1|,\dots, |k_{m_1}|\right\}\le A_1\log|n|+A_2.
\end{equation}
So, for each non-negative integer $n\le M$ (for some given upper bound M) for which there exist integers $k_i$ satisfying the equation~\eqref{eq:F sequence 10}, we have that $|k_i|\le A_2 + A_1\log M$, which means that we have at most $A_3\(1+\log M\)^{m_1}$ possible tuples $(k_1,\ldots, k_{m_1})\in\Z^{m_1}$, which may correspond to some $n\in\{0,\ldots, M\}$ solving the equation~\eqref{eq:F sequence 10} (where, once again, $A_3$ is a constant depending only on the initial data in our problem). Since $Q$ is a polynomial of degree $D\ge 1$, we conclude that the number of solutions $0\le n\le M$ to the equation~\eqref{eq:F sequence 10} is bounded above by $D\cdot A_3\left(1+\log M\right)^{m_1}$. Finally, recalling that $m_1\le m\le \dim V$, we obtain the desired conclusion from inequality~\eqref{eq:bound 3}.

{\bf Case 2.} $\mu_1\ne 1$.

In this case, since we also know that any characteristic root $\mu_r$ of the linear reccurence sequence $\{u_n\}_{n\in\bN_0}$ is either equal to $1$, or not a root of unity, we conclude that $\mu_1$ is not a root of unity.

The equation \eqref{eq:F sequence 3} reads now:
\begin{equation}
\label{eq:F sequence 11}
Q_1(n)\mu_1^n = \sum_{(i,j)\in \Sigma_2} c_{i,j}\lambda_j^{a_ik_i}.
\end{equation} 
We analyze again the subgroup $G_\Sigma\subseteq \Z^{1+m_1}$ containing the tuples $(f_0,f_1,\dots, f_{m_1})$ of integers satisfying the equations~\eqref{eq:same root}, i.e.,
\begin{equation}
\label{eq:same root 3}
\mu_1^{f_0}=\lambda_j^{a_if_i}\text{ for each }(i,j)\in\Sigma_2.
\end{equation} 
Because $\mu_1$ is not a root of unity and also each $\lambda_j$ is not a root of unity, while the $a_i$ are positive integers, we conclude that a nontrivial tuple $(f_0,f_1,\dots, f_{m_1})$ satisfying the equations~\eqref{eq:same root 3} must actually have each entry nonzero (i.e., $f_i\ne 0$ for each $i=0,\dots, m_1$). Therefore, each $\lambda_j^{a_i}$ is multiplicatively dependent with respect to $\mu_1$ and so, there exists an algebraic number $\lambda$ (which is not a root of unity), there exists a nonzero integer $b$ such that $\mu_1=\lambda^b$, and whenever there is a pair $(i,j)\in\Sigma_2$, there exist roots of unity $\zeta_{j,i}$ along with nonzero integers $b_i$ such that
\begin{equation}
\label{eq:same root 2}
\lambda_j^{a_i}=\zeta_{j,i}\cdot\lambda^{b_i}.
\end{equation}
We let $E$ be a positive integer such that $\zeta_{j,i}^E=1$ for each $(j,i)\in\Sigma_2$; then we let $B_i:=E\cdot b_i$ for each $i=1,\dots, m_1$. 
We now put  each exponent $k_i$ appearing in~\eqref{eq:F sequence 11}  in a prescribed  residue class modulo $E$ (just getting $E^m$ 
possible choices) and 
use~\eqref{eq:same root 2} along with the fact that $\mu_1=\lambda^b$. 
Writing $K_i: = \fl{k_i/E}$, $i=1, \ldots, m_1$, we obtain that 
 finding $n\in\bN_0$ which solves the 
equation~\eqref{eq:F sequence 11} (and then, in turn, also~\eqref{eq:F sequence 3} 
and~\eqref{eq:F sequence})  reduces to finding $n\in\bN_0$ which solves at least one of 
the at most $E^m$ 
distinct 
equations of the form:
\begin{equation}
\label{eq:equation 6}
 Q_1(n)\lambda^{bn}=\sum_{i=1}^{m_1} d_i\lambda^{B_i K_i},
\end{equation}
for some algebraic numbers $d_1,\ldots, d_{m_1}$, depending only on $E$, the $c_i$, and the $\zeta_{j,i}$, $(i,j) \in \Sigma_2$.  So, dividing the equation~\eqref{eq:equation 6} by $\lambda^{bn}$ yields that 
\begin{equation}
\label{eq:poly-expo} 
Q_1(n)=\sum_{i=1}^{m_1} d_i\lambda^{g_i},
\end{equation}
for some integers $g_i$.  Then once again applying \cite[Th\'{e}or\`{e}me~6, part~(b)]{Laurent} (see also our inequality~\eqref{eq:logarithmic growth}) yields immediately that any solution 
$(n,g_1,\ldots, g_{m_1})$ to the equation~\eqref{eq:poly-expo} must satisfy the inequality: 
$$
\max\{|g_1|,\dots, |g_{m_1}|\}\le A_4\log|n|+A_5,
$$
for some constants $A_4$ and $A_5$ depending only on the initial data in our problem ($X$, $\Phi$, $\alpha$, $V$). Then once again (exactly as in {\bf Case 1}), we conclude that there exists a constant $A_6$ such that for any  given upper bound $M\in\bN$, we have at most $A_6\(1+\log M\)^{m_1}$ possible tuples $(g_1,\ldots, g_{m_1})\in\Z^{m_1}$, which may correspond to some $n\in\{0,\ldots, M\}$ solving the equation~\eqref{eq:poly-expo}. Since $Q$ is a polynomial of degree $D\ge 1$, we conclude that the number of solutions $0\le n\le M$ to the equation~\eqref{eq:poly-expo} is bounded above by $D\cdot A_6\left(1+\log M\right)^{m_1}$. Finally, recalling that $m_1\le m\le \dim V$, we obtain the desired conclusion from inequality~\eqref{eq:bound 3}.

This concludes our proof of Theorem~\ref{thm:main 2}.

\section{Comments} 

\begin{remark}
If in the equation~\eqref{eq:F sequence} there exists at least one characteristic root $\mu_r$ of $\{u_n\}$ which is multiplicatively independent with respect to each one of the $\lambda_j$, then there is never a subsum~\eqref{eq:F sequence 3} containing $\mu_r$ on its left-hand side for which the corresponding group $G_\Sigma$ would be nontrivial. So, in this case, the equation~\eqref{eq:F sequence} would have only finitely many solutions. Therefore, with the notation as in Theorem~\ref{thm:main 2}, arguing as in the proof of \cite[Theorem~1.3]{CGSZ}, one concludes that if $\Phi$ is a group endomorphism of the semiabelian variety $X$ with the property that each characteristic root of its minimal polynomial (in ${\rm End}(X)$) is multiplicatively independent with respect to each eigenvalue $\lambda_j$ of the Frobenius endomorphism of $X$, then for each $\alpha\in X(K)$, the set $\mathcal{S}(\Phi, \alpha;V)$ defined by~\eqref{eq:DML} is a finite union of arithmetic progressions. 
\end{remark}

\begin{remark}
We notice that in~\eqref{eq:poly-expo}, if we deal with a polynomial $Q$ of degree $1$, then the conclusion from inequality~\eqref{eq:bound 3} is sharp. More precisely, as a specific example, the number of positive integers $n\le M$  which have precisely $m$ nonzero digits (all equal to $1$) in base-$p$ is 
of the order of $\(\log M\)^m$, which shows that Theorem~\ref{thm:main 2} is tight  if the Dynamical Mordell-Lang Conjecture reduces to solving the equation~\eqref{eq:poly-expo} when  $Q(n) = n$, 
$c_1=\ldots=c_m=1$ and $\lambda = p$. As proven in~\cite[Theorem~1.4]{CGSZ}, there are instances when the Dynamical Mordell-Lang Conjecture reduces \emph{precisely} to such equation.

Now, for higher degree polynomials $Q\in \Z[x]$ appearing in the equation~\eqref{eq:poly-expo}, one expects a lower exponent than $m$ appearing in the upper bounds from~\eqref{eq:bound 3}. One also notices that for any polynomial  $Q$, arguments $n$ with $k$ nonzero digits in base-$p$ lead to sparse outputs.  
Hence, simple combinatorics allows us to obtain a lower bound on the best 
possible exponent in~\eqref{eq:bound 3}.  
 However, finding a more precise exponent replacing $m$ in~\eqref{eq:bound 3} when $\deg Q>1$ seems very difficult beyond some special cases; the authors hope to return to this problem in a sequel paper.
\end{remark}

\section*{Acknowledgement}

 D.~G. and S~S. were partially supported by a Discovery Grant from NSERC, 
A.~O. by ARC Grants~DP180100201 and DP200100355, and I.~S.  by an ARC Grant~DP200100355.



\begin{thebibliography}{CGSZ20}


 
\bibitem[BGT15]{BGT-Funct}
J.~P.~Bell, D.~Ghioca, and T.~J.~Tucker, \emph{The Dynamical Mordell-Lang problem for Noetherian spaces}, Funct. Approx. Comment. Math. \textbf{53} (2015),   313--328.


\bibitem[BGT16]{BGT}
J.~P.~Bell, D.~Ghioca, and T.~J.~Tucker, \emph{The dynamical Mordell-Lang conjecture}, Mathematical Surveys and Monographs, \textbf{210}, American Mathematical Society, Providence, RI, 2016, xiii+280 pp.

\bibitem[CGSZ20]{CGSZ}
P.~Corvaja, D.~Ghioca, T.~Scanlon, and U.~Zannier, \emph{The Dynamical Mordell-Lang Conjecture for endomorphisms of semiabelian varieties defined over fields of positive characteristic}, J. Inst.  Math. Jussieu,  to appear. 

\bibitem[CZ13]{CZ3}
P.~Corvaja and U.~Zannier, \emph{Finiteness of odd perfect powers with four nonzero binary digits}, Ann. Inst. Fourier (Grenoble) \textbf{63} (2013),   715--731.

\bibitem[Ghi19]{G-TAMS}
D.~Ghioca, \emph{The dynamical Mordell-Lang conjecture in positive characteristic}, Trans. Amer. Math. Soc. \textbf{371} (2019),   1151--1167.  






\bibitem[GT09]{GT-JNT}
D.~Ghioca and T.~J.~Tucker, \emph{Periodic points, linearizing maps, and the dynamical Mordell-Lang problem}, J. Number Theory \textbf{129} (2009), 1392--1403.


\bibitem[Lau84]{Laurent}
M.~Laurent, \emph{\'Equations diophantiennes exponentielles}, Invent. Math. \textbf{78} (1984), 299--327.

\bibitem[MS04]{Moosa-Scanlon}
R.~Moosa and T.~Scanlon, \emph{$F$-structures and integral points on semiabelian varieties over finite fields}, Amer. J. Math. \textbf{126} (2004),  473--522.

\bibitem[Sch03]{Sch03}
W.~Schmidt, \emph{Linear recurrence sequences}, 
Diophantine Approximation (Cetraro, Italy, 2000),
Lecture Notes in Math. 1819, Springer-Verlag 
Berlin Heidelberg, 2003, pp.~171--247.


\end{thebibliography}
\end{document}